\documentclass{amsart}
\usepackage{cases}
\usepackage{mathrsfs}
\usepackage{graphicx}

\newtheorem{thm}{Theorem}[section]

\newtheorem{lem}[thm]{Lemma}

\theoremstyle{definition}
\newtheorem{defn}[thm]{Definition}
\theoremstyle{remark}
\newtheorem{rem}[thm]{Remark}
\numberwithin{equation}{section}

\begin{document}
\title[Multilinear commutators on weighted Morrey spaces]
{Estimates for multilinear commutators of generalized fractional integral operators on weighted Morrey spaces$^*$}
\author{Sha He and Xiangxing Tao}
\date{}






\address{school of mathematical sciences \endgraf
Beijing Normal University \endgraf Laboratory of mathematics and complex systems, ministry of education \endgraf Beijing, 100875,
P.R. China } \email{shahe@mail.bnu.edu.cn}

\address{department of mathematics \endgraf
Zhejiang University of Science and Technology \endgraf Hangzhou,
Zhejiang 310023, P.R. China } \email{xxtao@zust.edu.cn}

\begin{abstract}
Let $L$ be the infinitesimal generator of an analytic semigroup on
$L^2(\mathbb{R}^n)$ with Gaussian kernel bounds, and let $L^{-\alpha/2}$
be the fractional integrals of $L$ for $0<\alpha<n$. Assume that $\vec{b}=(b_1,b_2,\cdots,b_m)$ is a finite family of locally integrable functions, then the multilinear commutators generated by $\vec{b}$ and  $L^{-\alpha/2}$ is defined by
\begin{equation*}
L_{\vec{b}}^{-\alpha/2}f=[b_m,\cdots,[b_2,[b_1,L^{-\alpha/2}]],\cdots,]f
\end{equation*}
when $b_j\in BMO(w)$, $j=1,2,\cdots,m$, the authors obtain the boundedness of $L_{\vec{b}}^{-\alpha/2}$ on weighted Morrey spaces.
\end{abstract}

\maketitle

\section{Introduction and Main Results}

Assume that $L$ is a linear operator on $L^2(\mathbb{R}^n)$, which generates an analytic semigroup $e^{-tL}$ with a kernel $p_t(x,y)$ satisfying a Gaussian upper bound, that is,
\begin{align}\label{eq1}
|p_t(x,y)|\leq \frac{C}{t^{n/2}}e^{-c\frac{|x-y|^2}{t}},
\end{align}
for $x,y\in \mathbb{R}^n$ and all $t>0$.

For $0<\alpha<n$, the fractional integral $L^{-\alpha/2}$ generated by the operator $L$ is defined by
\begin{align*}
L^{-\alpha/2}f(x)=\frac{1}{\Gamma(\alpha/2)}\int_0^{\infty}e^{-tL}(f)\frac{dt}{t^{-\alpha/2+1}}(x).
\end{align*}

Note that if $L=-\Delta$, which is the Laplacian on $\mathbb{R}^n$, then $L^{-\alpha/2}$ is the classical fractional integral $I_{\alpha}$:
\begin{align*}
I_{\alpha}f(x)=\frac{\Gamma((n-\alpha)/2)}{\pi^{n/2}2^{\alpha}\Gamma(\alpha/2)}\int_{\mathbb{R}^n}\frac{f(y)}{|x-y|^{n-\alpha}}dy.
\end{align*}

Let $\vec{b}=(b_1,b_2,\cdots,b_m)$ be a finite family of locally integrable functions, then the multilinear commutators generated by $L^{-\alpha/2}$ and $\vec{b}$ is defined by
\begin{align*}
L_{\vec{b}}^{-\alpha/2}f=[b_m,\cdots,b_2,[b_1,L^{-\alpha/2}]],\cdots,]f
\end{align*}
where $m\in \mathbb{Z}^+$.

It is obvious that when $m=1$, denote $b_1=b$, then $L_{\vec{b}}^{-\alpha/2}f=[b,L^{-\alpha/2}]f$ is the commutator generated by $L^{-\alpha/2}$ and $b$, and when $b_1=b_2=\cdots=b_m$, $L_{\vec{b}}^{-\alpha/2}$ is the higher commutator.

It is well known that if $b\in BMO$, the commutator of fractional integral operator $[b,I_{\alpha}]$ is bounded from $L^p(\mathbb{R}^n)$ to $L^q(\mathbb{R}^n)$, where $1<p<n/\alpha$, $1/q=1/p-\alpha/n$, see \cite{Ch}. In 2004, Duong and Yan \cite{DY} proved the $(L^p,L^q)$ boundedness of commutator $[b,L^{-\alpha/2}]$ under the above conditions. Recently, Wang \cite{W2} obtained if $b\in\dot{\Lambda}_{\beta}(\mathbb{R}^n)$ (Lipschitz space), $[b,L^{-\alpha/2}]$ is of $(p,s)$ type, where $1/s=1/p-(\alpha+\beta)/n$; he also proved $[b,L^{-\alpha/2}]$ is bounded on weighted Lebesgue space when $b\in BMO$ and $b\in\dot{\Lambda}_{\beta}$. In \cite{ML}, Mo and Lu studied the $(L^p,L^q)$ boundedness of the operator $L_{\vec{b}}^{-\alpha/2}$ if $b_j\in BMO$ or $b_j\in \dot{\Lambda}_{\beta_j}$, $j=1,2,\cdots,m$.

On the other hand, the classical Morrey space was introduced by C.B.Morrey in the 1930s. Then there are lots of study of classical operators on Morrey spaces and Morrey type spaces, for detail see \cite{CF}, \cite{HT1}, \cite{HT2}, \cite{ST}, \cite{TS1}-\cite{TZheng}. In 2009, Komori and Shirai gave the definition of weighted Morrey space, and studied the boundedness of Hardy-Littlewood operator, fractional integral operator, Calder\'{o}n-Zygmund operator and its commutators on this space. After that, the study about this space has been increased. Wang \cite{W1} proved when $b\in BMO_1(w)$ (weighted BMO space) or $b\in\dot{\Lambda}_{\beta}(w)$ (weighted Lipschitz space), commutator $[b,I_{\alpha}]$ is bounded on weighted Morrey space. Soon, Si and Zhao \cite{S}, \cite{SZ} obtained the equivalent conditions of the boundedness of $[b,L^{-\alpha/2}]$ on weighted Morrey spaces. Therefore, they characterize weighted BMO space or weighted Lipschitz space by the boundedness of $[b,L^{-\alpha/2}]$ on this space, and generalize the results in \cite{W1}, \cite{W2}.

Motivated by the above results, in this paper, we will consider the local integrable function $b$ belongs to another weighted BMO space, that is, when $b_j\in BMO(w)$, $j=1,2,\cdots,m$, whose definition will be given later, we study the boundedness of $L_{\vec{b}}^{-\alpha/2}$ on weighted Morrey space. The main result of this paper is as follows.
\begin{thm}\label{thm1}
Assume the condition (\ref{eq1}) holds. Let $0<\alpha<n$, $1<p<n/\alpha$, $1/q=1/p-\alpha/n$, $0\leq \kappa <p/q$, $w^{q/p}\in A_1$, and $r_w>\frac{1-\kappa}{p/q-\kappa}$, where $r_w$ denotes the critical index of $w$ for the reverse H\"{o}lder condition. If $b_j\in BMO(w)$, $j=1,2,\cdots,m$, then
\begin{align*}
\|L_{\vec{b}}^{-\alpha/2}f\|_{L^{q,\kappa q/p}(w^{q/p},w)}\leq C\|\vec{b}\|_{BMO(w)}\|f\|_{L^{p,\kappa}(w)}.
\end{align*}
\end{thm}

\section{Definitions and lemmas}
A weight is a locally integrable function on $\mathbb{R}^n$ which takes values in $(0,\infty)$ almost everywhere. For a weight $w$ and a measurable set $E$, we define $w(E)=\int_E w(x)dx$, the Lebesgue measure of $E$ by $|E|$ and the characteristic function of $E$ by $\chi_E$. For a real number $p$, $1<p<\infty$, $p'$ is the conjugate of $p$, i.e. $1/p+1/p'=1$. The letter $C$ denotes a positive constant that may vary at each occurrence but is independent of the essential variable.

\begin{defn}
Let $1\leq p<\infty$, $0<\kappa<1$, $w$ is a weight, then weighted Morrey space is defined by
\begin{align*}
L^{p,\kappa}(w):=\{f\in L_{loc}^{p}(w):\|f\|_{L^{p,\kappa}(w)}<\infty\}
\end{align*}
where
\begin{equation*}
\|f\|_{L^{p,\kappa}(w)}=\sup\limits_{B}\Big(\frac{1}{w(B)^{\kappa}}\int_{B}|f(x)|^pw(x)dx\Big)^{\frac{1}{p}}
\end{equation*}
and the supremum is taken over all balls $B$ in $\mathbb{R}^n$.
\end{defn}

\begin{defn}
Let $1\leq p<\infty$, $0<\kappa<1$, $u,v$ are weight, then two weights weighted Morrey space is defined by
\begin{align*}
L^{p,\kappa}(u,v):=\{f:\|f\|_{L^{p,\kappa}(u,v)}<\infty\}
\end{align*}
where
\begin{equation*}
\|f\|_{L^{p,\kappa}(u,v)}=\sup\limits_{B}\Big(\frac{1}{v(B)^{\kappa}}\int_{B}|f(x)|^pu(x)dx\Big)^{\frac{1}{p}}
\end{equation*}
and the supremum is taken over all balls $B$ in $\mathbb{R}^n$. If $u=v$, then we denote $L^{p,\kappa}(u)$ for short.
\end{defn}

\begin{defn}
A weight function $w$ is in the Muckenhoupt class $A_p$ with $1<p<\infty$ if for every ball $B$ in $\mathbb{R}^n$, there exists a positive constant $C$ which is independent of $B$ such that
\begin{align*}
\Big(\frac{1}{|B|}\int_{B}w(x)dx\Big)\Big(\frac{1}{|B|}\int_{B}w(x)^{-\frac{1}{p-1}}dx\Big)^{p-1}\leq C
\end{align*}
When $p=1$, $w\in A_1$, if
\begin{align*}
\frac{1}{|B|}\int_{B}w(x)dx\leq C ess\inf\limits_{x\in B}w(x)
\end{align*}
When $p=\infty$, $w\in A_{\infty}$, if there exist positive constant $\delta$ and $C$ such that given a ball $B$ and $E$ is a measurable subset of $B$, then
\begin{align*}
\frac{w(E)}{w(B)}\leq \Big(\frac{|E|}{|B|}\Big)^{\delta}.
\end{align*}
\end{defn}

\begin{defn}
A weight function $w$ belongs to the reverse H\"{o}lder class $RH_r$ if there exist two constants $r>1$ and $C>0$ such that the following reverse H\"{o}lder inequality
\begin{align*}
\Big(\frac{1}{|B|}\int_Bw(x)^rdx\Big)^{\frac{1}{r}}\leq C\frac{1}{|B|}\int_Bw(x)dx.
\end{align*}
holds for every ball $B$ in $\mathbb{R}^n$.
\end{defn}
It is well known that if $w\in A_p$ with $1\leq p<\infty$, then there exists $r>1$ such that $w\in RH_r$. It follows from H\"{o}lder inequality that $w\in RH_r$ implies $w\in RH_s$ for all $1<s<r$. Moreover, if $w\in RH_r$, $r>1$, then we have $w\in RH_{r+\varepsilon}$ for some $\varepsilon>0$. We thus write $r_w\equiv\sup\{r>1:w\in RH_r\}$ to denote the critical index of $w$ for the reverse H\"{o}lder condition.

\begin{defn}
The Hardy-Littlewood maximal opoerator $M$ is defined by
\begin{align*}
Mf(x)=\sup\limits_{x\in B}\frac{1}{|B|}\int_B|f(y)|dy
\end{align*}
Let $w$ be a weight. The weighted maximal operator $M_w$ is defined by
\begin{align*}
M_wf(x)=\sup\limits_{x\in B}\frac{1}{w(B)}\int_B|f(y)|w(y)dy
\end{align*}
For $0<\alpha<n$, $r\geq 1$, the fractional maximal operator $M_{\alpha,r}$ is defined by
\begin{align*}
M_{\alpha,r}f(x)=\sup\limits_{x\in B}\Big(\frac{1}{|B|^{1-\alpha r/n}}\int_B|f(y)|^rdy\Big)^{1/r}
\end{align*}
and the fractional weighted maximal operator $M_{\alpha,r,w}$ is defined by
\begin{align*}
M_{\alpha,r,w}f(x)=\sup\limits_{x\in B}\Big(\frac{1}{w(B)^{1-\alpha r/n}}\int_B|f(y)|^rw(y)dy\Big)^{1/r}
\end{align*}
If $\alpha=0$, we denote $M_{r,w}$ for short.
\end{defn}

\begin{defn}
A family of operators $\{A_t:t>0\}$ is said to be an "approximation to identity" if, for every $t>0$, $A_t$ is represented by the kernel $a_t(x,y)$, which is a measurable function defined on $\mathbb{R}^n\times\mathbb{R}^n$, in the following sense: for every $f\in L^p(\mathbb{R}^n)$, $p\geq 1$,
\begin{align*}
A_tf(x)=\int_{\mathbb{R}^n}a_t(x,y)f(y)dy,
\end{align*}
and
\begin{align*}
|a_t(x,y)|\leq h_t(x,y)=t^{-n/2}g(\frac{|x-y|^2}{t}),
\end{align*}
for $(x,y)\in \mathbb{R}^n\times\mathbb{R}^n$, $t>0$. Here $g$ is a positive, bounded, decreasing function satisfying
\begin{align*}
\lim\limits_{r\rightarrow\infty}r^{n+\varepsilon}g(r^2)=0
\end{align*}
for some $\varepsilon>0$.
\end{defn}
Associated with an "approximation to identity" $\{A_t:t>0\}$, Martell \cite{M} introduced the sharp maximal function as follows:
\begin{align*}
M_A^{\sharp}f(x)=\sup\limits_{x\in B}\frac{1}{|B|}\int_B|f(y)-A_{t_B}f(y)|dy
\end{align*}
where $t_B=r_B^2$, $r_B$ is the radius of the ball $B$ and $f\in L^p(\mathbb{R}^n)$ for some $p\geq 1$.

Notice that our analytic semigroup $\{e^{-tL}:t>0\}$ is an "approximation to identity". In particular, denote
\begin{align*}
M_L^{\sharp}f(x)=\sup\limits_{x\in B}\frac{1}{|B|}\int_B|f(y)-e^{-t_BL}f(y)|dy.
\end{align*}

\begin{defn}\cite{MW}\label{2.7}
Let $w\in A_{\infty}$. Then the norm of $BMO(w)$: $\|\cdot\|_{*,w}$, is equivalent to the norm of $BMO(\mathbb{R}^n)$: $\|\cdot\|_{*}$, where
\begin{align*}
BMO(w)=\Big\{b:\|b\|_{*,w}=\sup\limits_{Q}\frac{1}{w(Q)}\int_Q|b(x)-b_{Q,w}|w(x)dx\Big\}
\end{align*}
and
\begin{align*}
b_{Q,w}=\frac{1}{w(Q)}\int_Qb(z)w(z)dz
\end{align*}
\end{defn}

Then, let us make some notations.

Given any positive integer $m$, for all $0\leq j\leq m$, we denote by $C_j^m$ the family of all finite subsets $\sigma=\{\sigma_1,\sigma_2,\cdots,\sigma_j\}$ of $\{1,2,\cdots,m\}$ of different elements, and for any $\sigma\in C_j^m$, let $\sigma'=\{1,2,\cdots,m\}\setminus\sigma$. Let $\vec{b}=(b_1,b_2,\cdots,b_m)$, then for any $\sigma=\{\sigma_1,\sigma_2,\cdots,\sigma_j\}\in C_j^m$, we denote $\vec{b}_{\sigma}=(b_{\sigma_1},b_{\sigma_2},\cdots,b_{\sigma_j})$, $b_{\sigma}(x)=\prod\limits_{\sigma_j\in\sigma}b_{\sigma_j}(x)$ and $\|\vec{b}_{\sigma}\|_{BMO(w)}=\prod\limits_{\sigma_j\in\sigma}\|b_{\sigma_j}\|_{BMO(w)}$, $\|\vec{b}\|_{BMO(w)}=\prod\limits_{\sigma_j\in\{1,2,\cdots,m\}}\|b_{\sigma_j}\|_{BMO(w)}
=\prod\limits_{j=1}^m\|b_j\|_{BMO(w)}$.

\begin{lem}\label{lem7}\cite{JN}
Let $s>1$, $1\leq p<\infty$ and $A_p^s=\{w:w^s\in A_p\}$. Then
\begin{align*}
A_p^s=A_{1+(p-1)/s}\bigcap RH_s.
\end{align*}
In particular, $A_1^s=A_1\bigcap RH_s$.
\end{lem}

\begin{lem}\label{lem8}\cite{M}
Assume that the semigroup $e^{-tL}$ has a kernel $p_t(x,y)$ which satisfies the upper bound (\ref{eq1}). Take $\lambda>0$, $f\in L_0^1(\mathbb{R}^n)$ and a ball $B_0$ such that there exists $x_0\in B_0$ with $Mf(x_0)\leq \lambda$. Then, for every $w\in A_{\infty}$, $0<\eta<1$, we can find $\gamma>0$ (independent of $\lambda, B_0, f, x_0$) and constant $C, r>0$ (which only depend on $w$)
\begin{align*}
w(\{x\in B_0:Mf(x)>A\lambda, M_L^{\sharp}f(x)\leq \gamma\lambda\})\leq C\eta^rw(B_0),
\end{align*}
Where $A>1$ is a fixed constant which depends only on $n$.
\end{lem}

As a result, by using the standard arguments, we have the following estimates:

For every $f\in L^{p,\kappa}(u,v)$, $1<p<\infty$, $0\leq \kappa<1$, if $u,v\in A_{\infty}$, then
\begin{equation*}
\|f\|_{L^{p,\kappa}(u,v)}\leq \|Mf\|_{L^{p,\kappa}(u,v)}\leq \|M_L^{\sharp}f\|_{L^{p,\kappa}(u,v)}
\end{equation*}

In particular, when $u=v=w$, $w\in A_{\infty}$, we have
\begin{equation*}
\|f\|_{L^{p,\kappa}(w)}\leq \|Mf\|_{L^{p,\kappa}(w)}\leq \|M_L^{\sharp}f\|_{L^{p,\kappa}(w)}
\end{equation*}

\begin{lem}\label{lem9}\cite{W1}
Let $0<\alpha<n$, $1<p<n/\alpha$, $1/q=1/p-\alpha/n$ and $w^{q/p}\in A_1$, if $0<\kappa<p/q$, $r_w>\frac{1-\kappa}{p/q-\kappa}$, then
\begin{align*}
\|M_{\alpha,1}f\|_{L^{q,\kappa q/p}(w^{q/p},w)}\leq C\|f\|_{L^{p,\kappa}(w)}
\end{align*}
It is also holds for $I_{\alpha}$.
\end{lem}

\begin{lem}\label{lem10}\cite{W1}
Let $0<\alpha<n$, $1<p<n/\alpha$, $1/q=1/p-\alpha/n$ and $w^{q/p}\in A_1$, if $0<\kappa<p/q$, $1<r<p$, $r_w>\frac{1-\kappa}{p/q-\kappa}$, then
\begin{align*}
\|M_{r,w}f\|_{L^{q,\kappa q/p}(w^{q/p},w)}\leq C\|f\|_{L^{q,\kappa q/p}(w^{q/p},w)}
\end{align*}
\end{lem}

\begin{lem}\label{lem11}\cite{W1}
Let $0<\alpha<n$, $1<p<n/\alpha$, $1/q=1/p-\alpha/n$, $0<\kappa<p/q$, $w\in A_{\infty}$. Then for $1<r<p$,
\begin{align*}
\|M_{\alpha,r,w}f\|_{L^{q,\kappa q/p}(w)}\leq C\|f\|_{L^{p,\kappa}(w)}
\end{align*}
\end{lem}

\begin{lem}\label{lem12}
Let $0<\alpha<n$, $1<p<n/\alpha$, $1/q=1/p-\alpha/n$ and $w^{q/p}\in A_1$, if $0\leq\kappa<p/q$, $r_w>\frac{1-\kappa}{p/q-\kappa}$, then
\begin{align*}
\|L^{-\alpha/2}f\|_{L^{q,\kappa q/p}(w^{q/p},w)}\leq C\|f\|_{L^{p,\kappa}(w)}
\end{align*}
\end{lem}
{\bf{Proof}}\quad Since semigroup $e^{-tL}$ has a kernel $p_t(x,y)$ that satisfies the upper bound (\ref{eq1}), it is easy to see that for $x\in \mathbb{R}^n$, $L^{-\alpha/2}f(x)\leq I_{\alpha}(|f|)(x)$. From the boundedness of $I_{\alpha}$ on weighted Morrey space (see Lemma \ref{lem9}), we get
\begin{align*}
\|L^{-\alpha/2}f\|_{L^{q,\kappa q/p}(w^{q/p},w)}\leq \|I_{\alpha}f\|_{L^{q,\kappa q/p}(w^{q/p},w)}\leq C\|f\|_{L^{p,\kappa}(w)}
\end{align*}

\begin{rem}\label{rem14}
Since $I_{\alpha}$ is of weak $(1,n/(n-\alpha))$ type, from the above proof, we can obtain $L^{-\alpha/2}$ is of weak $(1,n/(n-\alpha))$ type.
\end{rem}

\begin{lem}\label{lem14}\cite{DY}
Assume the semigroup $e^{-tL}$ has a kernel $p_t(x,y)$ which satisfies the upper bound (\ref{eq1}). Then for $0<\alpha<n$, the differential operator $L^{-\alpha/2}-e^{-tL}L^{-\alpha/2}$ has an associated kernel $\tilde{K}_{\alpha,t}(x,y)$ which satisfies
\begin{align*}
\tilde{K}_{\alpha,t}(x,y)\leq \frac{C}{|x-y|^{n-\alpha}}\frac{t}{|x-y|^2}.
\end{align*}
\end{lem}

\begin{lem}\label{lem15}
Assume the semigroup $e^{-tL}$ has a kernel $p_t(x,y)$ which satisfies the upper bound (\ref{eq1}), $b\in BMO(w)$, $w\in A_1$. Then for $f\in L^p(\mathbb{R}^n)$, $p>1$, $\sigma\in C_j^m(j=1,2,\cdots,m)$, $1<\tau<\infty$,
\begin{align*}
\sup\limits_{x\in B}\frac{1}{|B|}\int_B|e^{-t_BL}((b-b_B)_{\sigma}f)(y)|dy\leq C\|\vec{b}_{\sigma}\|_{BMO(w)}M_{\tau,w}f(x)
\end{align*}
where $t_B=r_B^2$, $r_B$ is the radius of $B$.
\end{lem}
{\bf{Proof}}\quad For any $f\in L^p(\mathbb{R}^n)$, $x\in B$, we have
\begin{align*}
&\frac{1}{|B|}\int_B|e^{-t_BL}((b-b_B)_{\sigma}f)(y)|dy\\
&\leq \frac{1}{|B|}\int_B\int_{\mathbb{R}^n}|p_{t_B}(y,z)||(b(z)-b_B)_{\sigma}f(z)|dzdy\\
&\leq \frac{1}{|B|}\int_B\int_{2B}|p_{t_B}(y,z)||(b(z)-b_B)_{\sigma}f(z)|dzdy\\
&\quad+ \frac{1}{|B|}\int_B\sum\limits_{k=1}^{\infty}\int_{2^{k+1}B\setminus 2^kB}|p_{t_B}(y,z)||(b(z)-b_B)_{\sigma}f(z)|dzdy:=M+N.
\end{align*}
Since for any $y\in B$, $z\in 2B$, from (\ref{eq1}), we get
\begin{align*}
|p_{t_B}(y,z)|\leq Ct_B^{-n/2}\leq C|2B|^{-1}
\end{align*}
Thus,
\begin{align*}
&M\leq \frac{C}{|2B|}\int_{2B}|(b(z)-b_B)_{\sigma}||f(z)|dz\\
&=\frac{C}{|2B|}\int_{2B}\prod\limits_{\sigma_j\in \sigma}|b_{\sigma_j}(z)-(b_{\sigma_j})_B||f(z)|dz
\end{align*}
For simplicity, we only consider the case of $j=2$. The case of $j>2$ is the same. Then
\begin{align*}
&M\leq \frac{C}{|2B|}\int_{2B}\big(|b_{\sigma_1}(z)-(b_{\sigma_1})_{2B}|+|(b_{\sigma_1})_{2B}-(b_{\sigma_1})_B|\big)\\
&\ \ \ \ \ \ \ \ \big(|b_{\sigma_2}(z)-(b_{\sigma_2})_{2B}|+|(b_{\sigma_2})_{2B}-(b_{\sigma_2})_B|\big)|f(z)|dz\\
&\leq \frac{C}{|2B|}\int_{2B}|b_{\sigma_1}(z)-(b_{\sigma_1})_{2B}||b_{\sigma_2}(z)-(b_{\sigma_2})_{2B}||f(z)|dz\\
&+\frac{C}{|2B|}\int_{2B}|b_{\sigma_1}(z)-(b_{\sigma_1})_{2B}||(b_{\sigma_2})_{2B}-(b_{\sigma_2})_B||f(z)|dz\\
&+\frac{C}{|2B|}\int_{2B}|(b_{\sigma_1})_{2B}-(b_{\sigma_1})_B||b_{\sigma_2}(z)-(b_{\sigma_2})_{2B}||f(z)|dz\\
&+\frac{C}{|2B|}\int_{2B}|(b_{\sigma_1})_{2B}-(b_{\sigma_1})_B||(b_{\sigma_2})_{2B}-(b_{\sigma_2})_B||f(z)|dz\\
&:=M_1+M_2+M_3+M_4.
\end{align*}

We split $M_1$ as follows.
\begin{align*}
&M_1\leq \frac{C}{|2B|}\int_{2B}\Big\{|b_{\sigma_1}(z)-(b_{\sigma_1})_{2B,w}|+|(b_{\sigma_1})_{2B,w}-(b_{\sigma_1})_{2B}|\Big\}\\
&\ \ \ \ \ \ \ \ \ \ \ \ \ \ \ \ \ \ \ \ \Big\{|b_{\sigma_2}(z)-(b_{\sigma_2})_{2B,w}|+|(b_{\sigma_2})_{2B,w}-(b_{\sigma_2})_{2B}|\Big\}|f(z)|dz\\
& =\frac{C}{|2B|}\int_{2B}|b_{\sigma_1}(z)-(b_{\sigma_1})_{2B,w}||b_{\sigma_2}(z)-(b_{\sigma_2})_{2B,w}||f(z)|dz\\
&+\frac{C}{|2B|}\int_{2B}|b_{\sigma_1}(z)-(b_{\sigma_1})_{2B,w}||(b_{\sigma_2})_{2B,w}-(b_{\sigma_2})_{2B}||f(z)|dz\\
&+\frac{C}{|2B|}\int_{2B}|(b_{\sigma_1})_{2B,w}-(b_{\sigma_1})_{2B}||b_{\sigma_2}(z)-(b_{\sigma_2})_{2B,w}||f(z)|dz\\
&+\frac{C}{|2B|}\int_{2B}|(b_{\sigma_1})_{2B,w}-(b_{\sigma_1})_{2B}||(b_{\sigma_2})_{2B,w}-(b_{\sigma_2})_{2B}||f(z)|dz\\
&:=M_{11}+M_{12}+M_{13}+M_{14}
\end{align*}

Choose $\tau_1, \tau_2, \tau, s>1$ that satisfy $1/\tau_1+1/\tau_2+1/\tau+1/s=1$. Then from H\"{o}lder's inequality and $w\in A_1$, we have
\begin{align*}
&M_{11}=\frac{C}{|2B|}\int_{2B}|b_{\sigma_1}(z)-(b_{\sigma_1})_{2B,w}|w(z)^{\frac{1}{\tau_1}}|b_{\sigma_2}(z)-(b_{\sigma_2})_{2B,w}|w(z)^{\frac{1}{\tau_2}}|f(z)|w(z)^{\frac{1}{\tau}}w(z)^{-1+\frac{1}{s}}dz\\
&\leq \frac{C}{|2B|}\Big(\int_{2B}|b_{\sigma_1}(z)-(b_{\sigma_1})_{2B,w}|^{\tau_1}w(z)dz\Big)^{\frac{1}{\tau_1}}\Big(\int_{2B}|b_{\sigma_2}(z)-(b_{\sigma_2})_{2B,w}|^{\tau_2}w(z)dz\Big)^{\frac{1}{\tau_2}}\\
&\qquad\Big(\int_{2B}|f(z)|^{\tau}w(z)dz\Big)^{\frac{1}{\tau}}\Big(\int_{2B}w(z)^{-s+1}dz\Big)^{\frac{1}{s}}\\
&\leq \frac{C}{|2B|}\|\vec{b}_{\sigma}\|_{*,w}M_{\tau,w}f(x)w(2B)^{1-1/s}w(x)^{-1+1/s}|2B|^{1/s}\\
&= C\|\vec{b}_{\sigma}\|_{*}M_{\tau,w}f(x)\Big(\frac{w(2B)}{|2B|}\Big)^{1-1/s}w(x)^{-1+1/s}\\
&\leq C\|\vec{b}_{\sigma}\|_*M_{\tau,w}f(x)
\end{align*}

In the above inequalities, we use the fact that if $w\in A_1$, then $w\in A_{\infty}$. Thus the norm of $BMO(w)$ is equivalent to the norm of $BMO(\mathbb{R}^n)$ (see Definition \ref{2.7}).\\

For $M_{12}$, we first estimate the term contains $(b_{\sigma_2})_{2B}$. In fact, it follows from the John and Nirenberg lemma that there exist $C_1>0$ and $C_2>0$ such that for any ball $B$ and $\alpha>0$
\begin{align*}
|\{z\in 2B: |b_{\sigma_2}(z)-(b_{\sigma_2})_{2B}|>\alpha\}|\leq C_1 |2B|e^{-C_2\alpha/\|b_{\sigma_2}\|_*}
\end{align*}
since $b\in BMO(\mathbb{R}^n)$. Using the definition of $A_{\infty}$, we get
\begin{align*}
w(\{z\in 2B: |b_{\sigma_2}(z)-(b_{\sigma_2})_{2B}|>\alpha\})\leq Cw(2B)e^{-C_2\alpha\delta/\|b_{\sigma_2}\|_*}
\end{align*}
for some $\delta>0$. Hence this inequality implies that
\begin{align*}
&\int_{2B}|b_{\sigma_2}(z)-(b_{\sigma_2})_{2B}|w(z)dz=\int_0^{\infty}w(\{z\in 2B: |b_{\sigma_2}(z)-(b_{\sigma_2})_{2B}|>\alpha\})d\alpha\\
&\ \ \ \ \ \ \ \ \ \ \ \ \ \ \ \ \ \ \ \ \ \ \ \ \ \ \ \ \ \ \ \ \ \ \ \ \ \ \leq Cw(2B)\int_0^{\infty}e^{-C_2\alpha\delta/\|b_{\sigma_2}\|_*}d\alpha\\
&\ \ \ \ \ \ \ \ \ \ \ \ \ \ \ \ \ \ \ \ \ \ \ \ \ \ \ \ \ \ \ \ \ \ \ \ \ \ =Cw(2B)\|b_{\sigma_2}\|_*.
\end{align*}
Thus,
\begin{align*}
&|(b_{\sigma_2})_{2B,w}-(b_{\sigma_2})_{2B}|\leq \frac{1}{w(2B)}\int_{2B}|b_{\sigma_2}(z)-(b_{\sigma_2})_{2B}|w(z)dz\\
&\ \ \ \ \ \ \ \ \ \ \ \ \ \ \ \ \ \ \ \ \ \ \ \ \ \ \ \leq \frac{C}{w(2B)}w(2B)\|b_{\sigma_2}\|_*=C\|b_{\sigma_2}\|_*\tag{2.1}
\end{align*}
We now estimate $M_{12}$ as follows: select $u$, $v$ such that $\frac{1}{u}+\frac{1}{v}+\frac{1}{\tau}=1$, by H\"{o}lder inequality, $w\in A_1$ and (2.1), we have
\begin{align*}
&M_{12}\leq C\|b_{\sigma_2}\|_*\frac{1}{|2B|}\int_{2B}|b_{\sigma_1}(z)-(b_{\sigma_1})_{2B,w}||f(z)|dz\\
&=C\|b_{\sigma_2}\|_*\frac{1}{|2B|}\int_{2B}|b_{\sigma_1}(z)-(b_{\sigma_1})_{2B,w}|w(z)^{\frac{1}{u}}|f(z)|w(z)^{\frac{1}{\tau}}w(z)^{-1+\frac{1}{v}}dz\\
&\leq C\|b_{\sigma_2}\|_*\frac{1}{|2B|}\Big(\int_{2B}|b_{\sigma_1}(z)-(b_{\sigma_1})_{2B,w}|^uw(z)dz\Big)^{\frac{1}{u}}\Big(\int_{2B}|f(z)|^{\tau}w(z)dz\Big)^{\frac{1}{\tau}}\Big(\int_{2B}w(z)^{-v+1}dz\Big)^{\frac{1}{v}}\\
&\leq C\|\vec{b}_{\sigma}\|_*M_{\tau,w}f(x)\frac{1}{|2B|}w(2B)^{1/u+1/\tau}w(x)^{-1+1/v}|2B|^{1/v}\\
&=C\|\vec{b}_{\sigma}\|_*M_{\tau,w}f(x)\Big(\frac{w(2B)}{|2B|}\Big)^{1-1/v}w(x)^{-1+1/v}\\
&=C\|\vec{b}_{\sigma}\|_*M_{\tau,w}f(x)
\end{align*}
The estimate for $M_{13}$ is the same as $M_{12}$, so we have $M_{13}\leq C\|\vec{b}_{\sigma}\|_*M_{\tau,w}f(x)$.\\
For the above $\tau$, choose $\tau'>1$ such that $\frac{1}{\tau}+\frac{1}{\tau'}=1$, then by H\"{o}lder inequality combining the inequality (2.1) and $w\in A_1$, we get
\begin{align*}
&M_{14}\leq C\|\vec{b}_{\sigma}\|_*\frac{1}{|2B|}\int_{2B}|f(z)|dz\\
&=C\|\vec{b}_{\sigma}\|_*\frac{1}{|2B|}\int_{2B}|f(z)|w(z)^{1/\tau}w(z)^{1/\tau'-1}dz\\
&\leq C\|\vec{b}_{\sigma}\|_*\frac{1}{|2B|}\Big(\int_{2B}|f(z)|^{\tau}w(z)dz\Big)^{1/\tau}\Big(\int_{2B}w(z)^{1-\tau'}dz\Big)^{1/\tau'}\\
&= C\|\vec{b}_{\sigma}\|_*\frac{1}{|2B|}M_{\tau,w}f(x)w(2B)^{1/\tau}w(x)^{1/\tau'-1}|2B|^{1/\tau'}\\
&= C\|\vec{b}_{\sigma}\|_*M_{\tau,w}f(x)\Big(\frac{w(2B)}{|2B|}\Big)^{1/\tau}w(x)^{1/\tau'-1}\\
&\leq C\|\vec{b}_{\sigma}\|_*M_{\tau,w}f(x)
\end{align*}
Combining the above results, we immediately have $M_1\leq C\|\vec{b}_{\sigma}\|_*M_{\tau,w}f(x)$.\\
For the term $M_2$, using the fact that $|b_{2B}-b_B|\leq C\|b\|_*$ (see \cite{T}), we now get
\begin{align*}
&M_2\leq \|b_{\sigma_2}\|_*\frac{C}{|2B|}\int_{2B}\Big(|b_{\sigma_1}(z)-(b_{\sigma_1})_{2B,w}|+|(b_{\sigma_1})_{2B,w}-(b_{\sigma_1})_{2B}|\Big)|f(z)|dz\\
&\ \ \ \ = \frac{C\|b_{\sigma_2}\|_*}{|2B|}\int_{2B}|b_{\sigma_1}(z)-(b_{\sigma_1})_{2B,w}||f(z)|dz\\
&\ \ \ \ +\frac{C\|b_{\sigma_2}\|_*}{|2B|}\int_{2B}|(b_{\sigma_1})_{2B,w}-(b_{\sigma_1})_{2B}||f(z)|dz\\
&\ \ \ \ :=M_{21}+M_{22}
\end{align*}

The estimate of $M_{21}$ is the same as $M_{12}$, thus we have $M_{21}\leq C\|b_{\sigma}\|_*M_{\tau,w}f(x)$.\\
For $M_{22}$, using the same method of $M_{14}$, we get $M_{22}\leq C\|b_{\sigma}\|_*M_{\tau,w}f(x)$.\\
Therefore,
\begin{equation*}
M_2\leq C\|b_{\sigma}\|_*M_{\tau,w}f(x)
\end{equation*}
Similarly,
\begin{equation*}
M_3\leq C\|b_{\sigma}\|_*M_{\tau,w}f(x)
\end{equation*}
Using the same method as $M_{14}$ to estimate $M_4$, we get
\begin{equation*}
M_4\leq C\|b_{\sigma}\|_*M_{\tau,w}f(x)
\end{equation*}
Hence, $M\leq C\|b_{\sigma}\|_*M_{\tau,w}f(x)$.\\

On the other hand, for any $y\in B$, $z\in 2^{k+1}B\setminus 2^kB$, we get $|y-z|\geq 2^{k-1}r_B$, then
\begin{align*}
|p_{t_B}(y,z)|\leq C\frac{e^{-C2^{2(k-1)}}2^{(k+1)n}}{|2^{k+1}B|}
\end{align*}
Thus,
\begin{align*}
&N\leq C\sum\limits_{k=1}^{\infty}\frac{e^{-C2^{2(k-1)}}2^{(k+1)n}}{|2^{k+1}B|}\int_{2^{k+1}B}|(b(z)-b_B)_{\sigma}||f(z)|dz\\
&=C\sum\limits_{k=1}^{\infty}\frac{e^{-C2^{2(k-1)}}2^{(k+1)n}}{|2^{k+1}B|}\int_{2^{k+1}B}|b_{\sigma_1}(z)-(b_{\sigma_1})_B|\cdots|b_{\sigma_j}(z)-(b_{\sigma_j})_B||f(z)|dz
\end{align*}
As before, for simplicity, we consider the case of $j=2$,
\begin{align*}
&N\leq C\sum\limits_{k=1}^{\infty}\frac{e^{-C2^{2(k-1)}}2^{(k+1)n}}{|2^{k+1}B|}\int_{2^{k+1}B}\big(|b_{\sigma_1}(z)-(b_{\sigma_1})_{2^{k+1}B}|+|(b_{\sigma_1})_{2^{k+1}B}-(b_{\sigma_1})_B|\big)\\
&\qquad\ \ \ \ \ \ \ \ \ \ \ \ \ \ \ \ \ \ \ \ \ \ \ \ \ \ \ \ \ \ \ \ \big(|b_{\sigma_2}(z)-(b_{\sigma_2})_{2^{k+1}B}|+|(b_{\sigma_2})_{2^{k+1}B}-(b_{\sigma_2})_B|\big)|f(z)|dz\\
&=C\sum\limits_{k=1}^{\infty}\frac{e^{-C2^{2(k-1)}}2^{(k+1)n}}{|2^{k+1}B|}\int_{2^{k+1}B}|b_{\sigma_1}(z)-(b_{\sigma_1})_{2^{k+1}B}||b_{\sigma_2}(z)-(b_{\sigma_2})_{2^{k+1}B}||f(z)|dz\\
&+C\sum\limits_{k=1}^{\infty}\frac{e^{-C2^{2(k-1)}}2^{(k+1)n}}{|2^{k+1}B|}\int_{2^{k+1}B}|(b_{\sigma_1})_{2^{k+1}B}-(b_{\sigma_1})_B||b_{\sigma_2}(z)-(b_{\sigma_2})_{2^{k+1}B}||f(z)|dz\\
&+C\sum\limits_{k=1}^{\infty}\frac{e^{-C2^{2(k-1)}}2^{(k+1)n}}{|2^{k+1}B|}\int_{2^{k+1}B}|b_{\sigma_1}(z)-(b_{\sigma_1})_{2^{k+1}B}||(b_{\sigma_2})_{2^{k+1}B}-(b_{\sigma_2})_B||f(z)|dz\\
&+C\sum\limits_{k=1}^{\infty}\frac{e^{-C2^{2(k-1)}}2^{(k+1)n}}{|2^{k+1}B|}\int_{2^{k+1}B}|(b_{\sigma_1})_{2^{k+1}B}-(b_{\sigma_1})_B||(b_{\sigma_2})_{2^{k+1}B}-(b_{\sigma_2})_B||f(z)|dz\\
&:=N_1+N_2+N_3+N_4.
\end{align*}
For $N_1$, similar to that of $M_1$, we have
\begin{align*}
&N_1\leq C\|\vec{b}_{\sigma}\|_{*}M_{\tau,w}f(x)\sum\limits_{k=1}^{\infty}e^{-C2^{2(k-1)}}2^{(k+1)n}\\
&\leq C\|\vec{b}_{\sigma}\|_{*}M_{\tau,w}f(x)
\end{align*}
Using the fact that $|b_{2^{j+1}B}-b_B|\leq 2^n (j+1)\|b\|_*$ (see \cite{T}), the estimates for $N_2, N_3, N_4$ are similar to $M_2, M_3, M_4$, thus
\begin{align*}
N_2, N_3, N_4\leq C\|\vec{b}_{\sigma}\|_{*}M_{\tau,w}f(x)
\end{align*}
Therefore, $N\leq C\|\vec{b}_{\sigma}\|_{*}M_{\tau,w}f(x)$.\\
Consequently, the lemma has been proved.

\begin{lem}\label{lem16}
Let $0<\alpha<n$, $w\in A_1$, $b\in BMO(w)$, then for all $r>1$, $\tau>1$ and $x\in\mathbb{R}^n$, we have
\begin{align*}
&M_L^{\sharp}(L_{\vec{b}}^{-\alpha/2}f)(x)\leq C\Big\{\|\vec{b}\|_{*}M_{r,w}(L^{-\alpha/2}f)(x)\\
&\ \ \ \ \ \ \ \ \ \ \ \ \ \ \ \ \ \ \ \ \ \ \ \ \ \ \ \ \ +\sum\limits_{j=1}^{m-1}\sum\limits_{\sigma\in C_j^m}C_{j,m}\|\vec{b}_{\sigma}\|_{*}M_{\tau,w}(L_{\vec{b}_{\sigma'}}^{-\alpha/2}f)(x)\\
&\ \ \ \ \ \ \ \ \ \ \ \ \ \ \ \ \ \ \ \ \ \ \ \ \ \ \ \ \  +\|\vec{b}\|_{*}w(x)^{-\alpha/n}M_{\alpha,r,w}f(x)+\|\vec{b}\|_{*}M_{\alpha,1}f(x)\Big\}
\end{align*}
\end{lem}
{\bf{Proof}}\quad For any given $x\in\mathbb{R}^n$, take a ball $B=B(x_0,r_B)$ which contains $x$. For $f\in L^p(\mathbb{R}^n)$, let $f_1=f\chi_{2B}$, $f_2=f-f_1$. Denote $K_{\alpha}(x,y)$ by the kernel of $L^{-\alpha/2}$, $\vec{\lambda}=(\lambda_1,\lambda_2,\cdots,\lambda_m)$, where $\lambda_j\in\mathbb{R}^n$, $j=1,2,\cdots,m$. Then $L_{\vec{b}}^{-\alpha/2}f$ can be written as the following form
\begin{align*}
&L_{\vec{b}}^{-\alpha/2}f(y)=\int_{\mathbb{R}^n}\prod\limits_{j=1}^m\big(b_j(y)-b_j(z)\big)K_{\alpha}(y,z)f(z)dz\\
& \ \ \ \ \ \ \ \ \ \ \ \ \ \ =\int_{\mathbb{R}^n}\prod\limits_{j=1}^m\big((b_j(y)-\lambda_j)-(b_j(z)-\lambda_j)\big)K_{\alpha}(y,z)f(z)dz\\
&\ \ \ \ \ \ \ \ \ \ \ \ \ \ =\sum\limits_{i=0}^m\sum\limits_{\sigma\in C_i^m}(-1)^{m-i}(b(y)-\lambda)_{\sigma}\int_{\mathbb{R}^n}(b(z)-\lambda)_{\sigma'}K_{\alpha}(y,z)f(z)dz\tag{2.2}
\end{align*}
Now expanding $(b(z)-\lambda)_{\sigma'}$ as
\begin{align*}
(b(z)-\lambda)_{\sigma'}=\big((b(z)-b(y))+(b(y)-\lambda)\big)_{\sigma'}
\end{align*}
Then it is easy to see from (2.2) that
\begin{align*}
&L_{\vec{b}}^{-\alpha/2}f(y)=\prod\limits_{j=1}^m(b_j(y)-\lambda_j)L^{-\alpha/2}f(y)+\sum\limits_{j=1}^{m-1}\sum\limits_{\sigma\in C_j^m}C_{j,m}(b(y)-\lambda)_{\sigma}L_{\vec{b}_{\sigma'}}^{-\alpha/2}f(y)\\
&\ \ \ \ \ \ \ \ \ \ +(-1)^mL^{-\alpha/2}\Big(\prod\limits_{j=1}^m(b_j-\lambda_j)f_1\Big)(y)+(-1)^mL^{-\alpha/2}\Big(\prod\limits_{j=1}^m(b_j-\lambda_j)f_2\Big)(y)
\end{align*}
where $C_{j,m}$ is a constant only relevant to $j, m$.\\
Thus
\begin{align*}
&e^{-t_BL}(L_{\vec{b}}^{-\alpha/2}f)(y)=e^{-t_BL}\Big(\prod\limits_{j=1}^m(b_j-\lambda_j)L^{-\alpha/2}f\Big)(y)\\
&\ \ \ \ \ \ \ \ \ \ \ \ \ \ \ \ \ \ \ \ \ \ \quad +\sum\limits_{j=1}^{m-1}\sum\limits_{\sigma\in C_j^m}C_{j,m}e^{-t_BL}\big((b-\lambda)_{\sigma}L_{\vec{b}_{\sigma'}}^{-\alpha/2}f\big)(y)\\
&\ \ \ \ \ \ \ \ \ \ \ \ \ \ \ \ \ \ \ \ \ \ \quad +(-1)^me^{-t_BL}\Big(L^{-\alpha/2}\big(\prod\limits_{j=1}^m(b_j-\lambda_j)f_1\big)\Big)(y)\\
&\ \ \ \ \ \ \ \ \ \ \ \ \ \ \ \ \ \ \ \ \ \ \quad+
(-1)^me^{-t_BL}\Big(L^{-\alpha/2}\big(\prod\limits_{j=1}^m(b_j-\lambda_j)f_2\big)\Big)(y)\\
\end{align*}
where $t_B=r_B^2$, $r_B$ is the radius of ball $B$.

Take $\lambda_j=(b_j)_B$, $j=1,2,\cdots,m$, and denote $\vec{b}_B=((b_1)_B,(b_2)_B,\cdots,(b_m)_B)$, then
\begin{align*}
&\frac{1}{|B|}\int_B\Big|L_{\vec{b}}^{-\alpha/2}f(y)-e^{-t_BL}(L_{\vec{b}}^{-\alpha/2}f)(y)\Big|dy\\
&\leq \frac{1}{|B|}\int_B\Big|\prod\limits_{j=1}^m(b_j(y)-(b_j)_B)L^{-\alpha/2}f(y)\Big|dy\\
&+\frac{1}{|B|}\int_B\Big|\sum\limits_{j=1}^{m-1}\sum\limits_{\sigma\in C_j^m}C_{j,m}(b(y)-b_B)_{\sigma}L_{\vec{b}_{\sigma'}}^{-\alpha/2}f(y)\Big|dy\\
&+\frac{1}{|B|}\int_B\Big|L^{-\alpha/2}\Big(\prod\limits_{j=1}^m(b_j-(b_j)_B)f_1\Big)(y)\Big|dy\\
&+\frac{1}{|B|}\int_B\Big|e^{-t_BL}(\prod\limits_{j=1}^m(b_j-(b_j)_B)L^{-\alpha/2}f)(y)\Big|dy\\
&+\frac{1}{|B|}\int_B\Big|\sum\limits_{j=1}^{m-1}\sum\limits_{\sigma\in C_j^m}C_{j,m}e^{-t_BL}((b-b_B)_{\sigma}L_{\vec{b}_{\sigma'}}^{-\alpha/2}f)(y)\Big|dy\\
&+\frac{1}{|B|}\int_B\Big|e^{-t_BL}\Big(L^{\alpha/2}\big(\prod\limits_{j=1}^m(b_j-(b_j)_B)f_1\big)\Big)(y)\Big|dy\\
&+\frac{1}{|B|}\int_B\Big|L^{\alpha/2}\Big(\prod\limits_{j=1}^m(b_j-(b_j)_B)f_2\Big)(y)-e^{-t_BL}\Big(L^{-\alpha/2}\big(\prod\limits_{j=1}^m(b_j-(b_j)_B)f_2\big)\Big)(y)\Big|dy\\
&:=I+II+III+IV+V+VI+VII.
\end{align*}
Let us now estimate $I,II,III,IV,V,VI,VII$ respectively.\\

We take $m=2$ as an example, the estimate for the case $m>2$ is the same. We split $I$ as follows:
\begin{align*}
&I=\frac{1}{|B|}\int_B|b_1(y)-(b_1)_{B,w}||b_2(y)-(b_2)_{B,w}||L^{-\alpha/2}f(y)|dy\\
&\ \ +\frac{1}{|B|}\int_B|b_1(y)-(b_1)_{B,w}||(b_2)_{B,w}-(b_2)_B||L^{-\alpha/2}f(y)|dy\\
&\ \ +\frac{1}{|B|}\int_B|(b_1)_{B,w}-(b_1)_B||b_2(y)-(b_2)_{B,w}||L^{-\alpha/2}f(y)|dy\\
&\ \ +\frac{1}{|B|}\int_B|(b_1)_{B,w}-(b_1)_B||(b_2)_{B,w}-(b_2)_B||L^{-\alpha/2}f(y)|dy\\
&\ \ :=I_1+I_2+I_3+I_4
\end{align*}
Choose $r_1,r_2,r,q>1$, such that $1/r_1+1/r_2+1/r+1/q=1$, then by H\"{o}lder's inequality, $w\in A_1$, and using the same estimate of $M_1$, we have
\begin{align*}
I\leq C\|\vec{b}\|_{*}M_{r,w}(L^{-\alpha/2}f)(x).
\end{align*}

For $II$, take $\tau_1,\cdots,\tau_j,\tau,\nu>1$ that satisfy $1/\tau_1+\cdots+1/\tau_j+1/\tau+1/\nu=1$, then by the same estimate as $I$, we get
\begin{align*}
II\leq C\sum\limits_{j=1}^{m-1}\sum\limits_{\sigma\in C_j^m}C_{j,m}\|\vec{b}_{\sigma}\|_*M_{\tau,w}(L_{\vec{b}_{\sigma'}}^{-\alpha/2}f)(x)
\end{align*}

Before we estimate $III$, let us introduce Kolmogorov's inequality (\cite{GF},P455):\\
Let $0<r<l<\infty$, for $f\geq 0$, define $\|f\|_{L^{l,\infty}}=\sup\limits_{t>0}t|\{x\in\mathbb{R}^n:|f(x)|>t\}|^{\frac{1}{l}}$, $N_{l,r}(f)=\sup\limits_E\frac{\|f\chi_E\|_r}{\|\chi_E\|_h}$, $\frac{1}{h}=\frac{1}{r}-\frac{1}{l}$, then
\begin{align*}
\|f\|_{L^{l,\infty}}\leq N_{l,r}(f)\leq \Big(\frac{l}{l-r}\Big)^{\frac{1}{r}}\|f\|_{L^{l,\infty}}.
\end{align*}
Applying Kolmogorov's inequality, weak $(1,n/(n-\alpha))$ boundedness of $L^{-\alpha/2}$ (see Remark \ref{rem14}) and H\"{o}lder's inequality, we have
\begin{align*}
&III=\frac{1}{|B|}\int_B|L^{-\alpha/2}(\prod\limits_{j=1}^m(b_j-(b_j)_B)f_1)(y)|dy\\
&\ \ \ \ \leq \frac{C}{|B|^{1-\alpha/n}}\|L^{-\alpha/2}(\prod\limits_{j=1}^m(b_j-(b_j)_B)f_1)\|_{L^{\frac{n}{n-\alpha},\infty}}\\
&\ \ \ \ \leq \frac{C}{|B|^{1-\alpha/n}}\int_{2B}|\prod\limits_{j=1}^m(b_j(y)-(b_j)_B)f(y)|dy
\end{align*}
We consider the case of $m=2$ for example. Take $r_1,r_2,r,q>1$ such that $1/r_1+1/r_2+1/r+1/q=1$, then
\begin{align*}
&III\leq \frac{C}{|B|^{1-\alpha/n}}\int_{2B}(|b_1(y)-(b_1)_{2B}|+|(b_1)_{2B}-(b_1)_B|)(|b_2(y)-(b_2)_{2B}|+|(b_2)_{2B}-(b_2)_B|)|f(y)|dy\\
&\ \ \ \ \leq \frac{C}{|B|^{1-\alpha/n}}\int_{2B}|b_1(y)-(b_1)_{2B}||b_2(y)-(b_2)_{2B}||f(y)|dy\\
&\ \ \ \ +\frac{C}{|B|^{1-\alpha/n}}\int_{2B}|b_1(y)-(b_1)_{2B}||(b_2)_{2B}-(b_2)_B|)||f(y)|dy\\
&\ \ \ \ +\frac{C}{|B|^{1-\alpha/n}}\int_{2B}|(b_1)_{2B}-(b_1)_B||b_2(y)-(b_2)_{2B}||f(y)|dy\\
&\ \ \ \ +\frac{C}{|B|^{1-\alpha/n}}\int_{2B}|(b_1)_{2B}-(b_1)_B||(b_2)_{2B}-(b_2)_B||f(y)|dy\\
&\ \ \ \ :=III_1+III_2+III_3+III_4.
\end{align*}
The estimates of $III_1,III_2,III_3$ are similar to that of $M_1,M_2,M_3$, thus we get
\begin{align*}
III_1,III_2,III_3\leq C\|\vec{b}\|_{*}M_{\alpha,r,w}f(x)
\end{align*}
For $III_4$:
\begin{align*}
&III_4\leq C\|\vec{b}\|_{*}\frac{1}{|B|^{1-\alpha/n}}\int_{2B}|f(y)|dy\\
&\ \ \ \ \ \ \leq C\|\vec{b}\|_{*}M_{\alpha,1}f(x)
\end{align*}
Therefore,
\begin{align*}
III\leq C\|\vec{b}\|_{*}M_{\alpha,r,w}f(x)+C\|\vec{b}\|_{*}M_{\alpha,1}f(x).
\end{align*}
By Lemma \ref{lem15}:
\begin{align*}
IV\leq C\|\vec{b}\|_{*}M_{r,w}(L^{-\alpha/2}f)(x)
\end{align*}
and
\begin{align*}
V\leq C\sum\limits_{j=1}^{m-1}\sum\limits_{\sigma\in C_j^m}C_{j,m}\|\vec{b}_{\sigma}\|_{*}M_{\tau,w}(L_{\vec{b}_{\sigma'}}^{-\alpha/2}f)(x)
\end{align*}
We estimate $VI$ as follows:
\begin{align*}
&VI=\frac{1}{|B|}\int_B\Big|\int_{\mathbb{R}^n}p_{t_B}(y,z)\Big(L^{-\alpha/2}\big(\prod\limits_{j=1}^m(b_j-(b_j)_B)f_1\big)\Big)(z)dz\Big|dy\\
&\leq \frac{1}{|B|}\int_B\int_{2B}\Big|p_{t_B}(y,z)\Big|\Big|L^{-\alpha/2}\Big(\prod\limits_{j=1}^m(b_j-(b_j)_B)f_1\Big)(z)\Big|dzdy\\
&\leq \frac{1}{|B|}\int_B\int_{\mathbb{R}^n\setminus {2B}}\Big|p_{t_B}(y,z)\Big|\Big|L^{-\alpha/2}\Big(\prod\limits_{j=1}^m(b_j-(b_j)_B)f_1\Big)(z)\Big|dzdy\\
&=VI_1+VI_2
\end{align*}
For $VI_1$, since $y\in B$, $z\in 2B$, $|p_{t_B}(y,z)|\leq C|2B|^{-1}$, thus
\begin{align*}
VI_1\leq\frac{C}{|2B|}\int_{2B}|L^{-\alpha/2}(\prod\limits_{j=1}^m(b_j-(b_j)_B)f_1)(z)|dz
\end{align*}
Then the estimate is similar to that of $III$, thus,
\begin{align*}
VI_1\leq C\|\vec{b}\|_{*}M_{\alpha,r,w}f(x)+C\|\vec{b}\|_{*}M_{\alpha,1}f(x)
\end{align*}
For $VI_2$, since $y\in B$, $z\in 2^{k+1}B\setminus 2^kB$, $|p_{t_B}(y,z)|\leq C\frac{e^{-C2^{2(k-1)}}2^{(k+1)n}}{|2^{k+1}B|}$. Hence,
\begin{align*}
VI_2\leq \sum\limits_{k=1}^{\infty}\frac{e^{-C2^{2(k-1)}}2^{(k+1)n}}{|2^{k+1}B|}\int_{2^{k+1}B}|L^{-\alpha/2}(\prod\limits_{j=1}^m(b_j-(b_j)_B)f_1)(z)|dz
\end{align*}
and following the same method of $III$, we have
\begin{align*}
&VI_2\leq \sum\limits_{k=1}^{\infty}\frac{e^{-C2^{2(k-1)}}2^{(k+1)n}}{(2^{k+1})^{n-\alpha}}\Big(\|\vec{b}\|_{*}M_{\alpha,r,w}f(x)+\|\vec{b}\|_{*}M_{\alpha,1}f(x)\Big)\\
&\leq C\|\vec{b}\|_{*}M_{\alpha,r,w}f(x)+C\|\vec{b}\|_{*}M_{\alpha,1}f(x)
\end{align*}
Applying Lemma \ref{lem14},
\begin{align*}
&VII\leq\frac{1}{|B|}\int_B|(L^{-\alpha/2}-e^{-t_BL}L^{-\alpha/2})(\prod\limits_{j=1}^m(b_j-(b_j)_B)f_2)(y)|dy\\
&\leq\frac{1}{|B|}\int_B\int_{\mathbb{R}^n\setminus {2B}}|\tilde{K}_{\alpha,t_B}(y,z)||\prod\limits_{j=1}^m(b_j(z)-(b_j)_B)f(z)|dzdy\\
&\leq C\sum\limits_{k=1}^{\infty}\int_{2^{k+1}B\setminus 2^kB}\frac{t_B}{|x_0-z|^{n-\alpha+2}}|\prod\limits_{j=1}^m(b_j(z)-(b_j)_B)f(z)|dz\\
&\leq C\sum\limits_{k=1}^{\infty}\frac{2^{-2k}}{|2^kB|^{1-\alpha/n}}\int_{2^{k+1}B}|\prod\limits_{j=1}^m(b_j(z)-(b_j)_B)||f(z)|dz
\end{align*}
Using the similar estimate of $III$, we have
\begin{align*}
&VII\leq C\sum\limits_{k=1}^{\infty}2^{-2k}\Big(\|\vec{b}\|_{*}M_{\alpha,r,w}f(x)+\|\vec{b}\|_{*}M_{\alpha,1}f(x)\Big)\\
&\leq C\|\vec{b}\|_{*}\Big(M_{\alpha,r,w}f(x)+M_{\alpha,1}f(x)\Big)
\end{align*}
Therefore, we have completed the proof of Lemma \ref{lem16}.

\section{The proof of theorem \ref{thm1}}
We are now in the position of proving Theorem \ref{thm1}.\\

{\bf{Proof}}\quad From Lemma \ref{lem7}, Lemma \ref{lem8}, Lemma \ref{lem16} and Lemma \ref{lem9}-Lemma \ref{lem12}, we obtain
\begin{align*}
&\|L_{\vec{b}}^{-\alpha/2}f\|_{L^{q,\kappa q/p}(w^{q/p},w)}\\
&\leq \|M_L^{\sharp}(L_{\vec{b}}^{-\alpha/2}f)\|_{L^{q,\kappa q/p}(w^{q/p},w)}\\
&\leq C\Big\{\|\vec{b}\|_{*}\|M_{r,w}(L^{-\alpha/2}f)\|_{L^{q,\kappa q/p}(w^{q/p},w)}\\
&\ \ \ \ \ \ \ +\sum\limits_{j=1}^{m-1}\sum\limits_{\sigma\in C_j^m}C_{j,m}\|\vec{b}_{\sigma}\|_{*}\|M_{\tau,w}(L_{\vec{b}_{\sigma'}}^{-\alpha/2}f)\|_{L^{q,\kappa q/p}(w^{q/p},w)}\\
&\ \ \ \ \ \ \ +\|\vec{b}\|_{*}\|w^{-\alpha/n}M_{\alpha,r,w}(f)\|_{L^{q,\kappa q/p}(w^{q/p},w)}+\|\vec{b}\|_{*}\|M_{\alpha,1}(f)\|_{L^{q,\kappa q/p}(w^{q/p},w)}\Big\}\\
&\leq C\|\vec{b}\|_{*}\|L^{-\alpha/2}f\|_{L^{q,\kappa q/p}(w^{q/p},w)}\\
&\ \ \ \ \ \ \ +\sum\limits_{j=1}^{m-1}\sum\limits_{\sigma\in C_j^m}C_{j,m}\|\vec{b}_{\sigma}\|_{*}\|(L_{\vec{b}_{\sigma'}}^{-\alpha/2}f)\|_{L^{q,\kappa q/p}(w^{q/p},w)}\\
&\ \ \ \ \ \ \ +C\|\vec{b}\|_{*}\|M_{\alpha,r,w}f\|_{L^{q,\kappa q/p}(w)}+C\|\vec{b}\|_*\|f\|_{L^{p,\kappa}(w)}\\
&\leq C\|\vec{b}\|_{*}\|f\|_{L^{p,\kappa}(w)}+\sum\limits_{j=1}^{m-1}\sum\limits_{\sigma\in C_j^m}C_{j,m}\|\vec{b}_{\sigma}\|_{*}\|(L_{\vec{b}_{\sigma'}}^{-\alpha/2}f)\|_{L^{q,\kappa q/p}(w^{q/p},w)}
\end{align*}
Then, we can make use of induction on $\sigma\subseteq \{1,2,\cdots,m\}$ to get that
\begin{align*}
\|L_{\vec{b}}^{-\alpha/2}f\|_{L^{q,\kappa q/p}(w^{q/p},w)}\leq C\|\vec{b}\|_{*}\|f\|_{L^{p,\kappa}(w)}.
\end{align*}
This completes the proof of Theorem \ref{thm1}.

\bigskip\

\end{document}